\input amstex
\mag=\magstephalf 
\documentstyle{amsppt}
\nologo
\define\BaseLineSkip{\baselineskip=14pt} 
\NoBlackBoxes
\rightheadtext\nofrills{HYPERELLIPTIC FUNCTIONS IN GENUS THREE}
\leftheadtext\nofrills{YOSHIHIRO \^ONISHI}
\define\inbox#1{$\boxed{\text{#1}}$}
\def\fp{\flushpar}
\define\dint{\dsize\int}

\define\underbarl#1{\lower 1.4pt \hbox{\underbar{\raise 1.4pt \hbox{#1}}}}

\def\fp{\flushpar}
\define\tp#1{\negthinspace\ ^t#1}

\define\lr#1{^{\sssize\left(#1\right)}}

\define\nullhbox{\hbox{\vrule height0pt depth 0pt width 1pt}}
\define\qedright{\null\nobreak\leaders\nullhbox\hskip10pt plus1filll\ \qed
\par $ $ \par}

\font\sc=cmcsc10
\font\twelveptrm=cmr12

\font\sc=cmcsc10
\BaseLineSkip
\topmatter
\title  \nofrills\twelveptrm 
DETERMINANT EXPRESSIONS FOR HYPERELLIPTIC FUNCTIONS IN GENUS THREE
\endtitle
\author {YOSHIHIRO \^ONISHI} 
\endauthor
\subjclass  11G30, 11G10, 14H45 \endsubjclass
\endtopmatter
\document
\TagsOnRight
\BaseLineSkip
%
%
%

\heading {\rm 1. {\it Introduction}}  \endheading

\fp
Let  $\sigma(u)$  and  $\wp(u)$  be the usual functions 
in the theory of elliptic functions.  
The following two formulae were found in the nineteenth-century.  
First one is 
  $$
  \aligned
  (-1)^{n(n-1)/2}1!2!&\cdots n!
   \frac{\sigma(u_0+u_1+\cdots+u_n)\prod_{i<j}\sigma(u_i-u_j)}
        {\sigma(u_0)^{n+1}\sigma(u_1)^{n+1}\cdots\sigma(u_n)^{n+1}} \\
  &= \left|\matrix
      1 &  \wp(u_0) & \wp'(u_0) & \wp''(u_0) & \cdots & \wp\lr{n-1}(u_0) \\
      1 &  \wp(u_1) & \wp'(u_1) & \wp''(u_1) & \cdots & \wp\lr{n-1}(u_1) \\
 \vdots &  \vdots   & \vdots    & \vdots     & \ddots & \vdots           \\
      1 &  \wp(u_n) & \wp'(u_n) & \wp''(u_n) & \cdots & \wp\lr{n-1}(u_n) \\
   \endmatrix\right|.
  \endaligned
  \tag 1.1
  $$
This formula appeared in  
the paper of Frobenius and Stickelberger \cite{{\bf 7}}.  
Second one is 
  $$
  (-1)^{n(n-1)/2}(1!2!\cdots (n-1)!)^2
   \frac{\sigma(nu)}{\sigma(u)^{n^2}}
  = \left|\matrix
    \wp'        & \wp''     &  \cdots  & \wp\lr{n-1}  \\
    \wp''       & \wp'''    &  \cdots  & \wp\lr{n}    \\
    \vdots      & \vdots    &  \ddots  & \vdots       \\
    \wp\lr{n-1} & \wp\lr{n} &  \cdots  & \wp\lr{2n-3} \\
   \endmatrix\right|(u).   
  \tag 1.2
  $$
Although this formula can be obtained by a limiting process from (1.1), 
it was found before \cite{{\bf 7}} by the paper of Kiepert \cite{{\bf 9}}.  

If we set  $y(u)=\frac 12 \wp'(u)$  and  $x(u)=\wp(u)$,   
then we have an equation  $y(u)^2=x(u)^3+\cdots$, that is a defining equation 
of the elliptic curve 
to which the functions  $\wp(u)$  and  $\sigma(u)$  are attached.  
Here the complex number  $u$  and the coordinate  $(x(u), y(u))$  
correspond by the equality  
  $$
  u=\dint_{\infty}^{(x(u), y(u))}\dfrac{dx}{2y}.
  $$
Then (1.1) and (1.2) is easily rewitten as
  $$
  \aligned
  &\frac{\sigma(u_0+u_1+\cdots+u_n)\prod_{i<j}\sigma(u_i-u_j)}
        {\sigma(u_0)^{n+1}\sigma(u_1)^{n+1}\cdots\sigma(u_n)^{n+1}} \\
  &=\left|\matrix
     1 &  x(u_0) & y(u_0) & x^2(u_0) & yx(u_0) & x^3(u_0) & \cdots  \\
     1 &  x(u_1) & y(u_1) & x^2(u_1) & yx(u_1) & x^3(u_1) & \cdots  \\
\vdots &  \vdots & \vdots   & \vdots & \vdots   & \vdots  & \ddots  \\
     1 &  x(u_n) & y(u_n) & x^2(u_n) & yx(u_n) & x^3(u_n) & \cdots  \\
   \endmatrix\right|
  \endaligned
  \tag 1.3
  $$
and
  $$
  \align
  1!2!&\cdots (n-1)!
  \frac{\sigma(nu)}{\sigma(u)^{n^2}} \\
  &=\left|\matrix
    x'           & y'            & (x^2)'  
  & (yx)'        & (x^3)'        & \cdots  \\
    x''          & y''           & (x^2)'' 
  & (yx)''       & (x^3)''       & \cdots  \\
    \vdots       & \vdots        & \vdots   
  & \vdots       & \vdots        & \ddots  \\
    x\lr{n-1}    & y\lr{n-1}     & (x^2)\lr{n-1} 
  & (yx)\lr{n-1} & (x^3)\lr{n-1} & \cdots  \\
  \endmatrix\right|(u), 
  \tag 1.4
  \endalign
  $$
respectively.  

The author recently gave a generalization of the formulae (1.3) and (1.4) to 
the case of genus two in \cite{{\bf 13}}. 
Our aim is to give a quite natural genaralization of  (1.3)  and  (1.4) 
and the results in \cite{{\bf 13}}  to the case of genus three 
(see Theorem 3.2 and Theorem 4.2).  
Our generalization of the function in the left hand side of (1.4) 
is along a line which appeared for a curve of genus two 
in the paper \cite{{\bf 8}} of D. Grant.  
Although Fay's famous formula, that is (44) in p.33 of \cite{{\bf 6}},  
possibly relates with our generalizations, 
no connections are known.  

Now we prepare the minimal fundamentals to explain our results.  
Let  $f(x)$  be a monic polynomial of  $x$  of degree $7$.  
Assume that  $f(x)=0$  has no multiple roots.   
Let  $C$  be the hyperelliptic curve defined by  $y^2=f(x)$.  
Then  $C$  is of genus  $3$ and it is ramified at infinity.  
We denote by $\infty$  the unique point at infinity.  
Let  $\bold C^3$  be the coordinate space of all values of 
the integrals, with their initial points  $\infty$,  of the first kind 
with respect to the basis  $dx/2y$, $xdx/2y$, $x^2dx/2y$  
of the differentials of first kind.  
Let  $\Lambda\subset \bold C^3$  be the lattice of their periods.  
So  $\bold C^3/\Lambda$  is the Jacobian variety of  $C$.  
We have an embedding  $\iota:C\hookrightarrow \bold C^3/\Lambda$  
defined 
by  $P\mapsto (\int_{\infty}^P\frac{dx}{2y}, \ 
               \int_{\infty}^P\frac{xdx}{2y}, \  
               \int_{\infty}^P\frac{x^2dx}{2y})$.  
Therefore  $\iota(\infty)=(0, 0, 0)\in \bold C^3/\Lambda$.  
We also have 
a canonical map  $\kappa: \bold C^3 \rightarrow \bold C^3/\Lambda$.  
An algebraic function on  $C$,  
that we call a {\it hyperelliptic function} in this article,  
is regarded as a function on a universal covering  $\kappa^{-1}\iota(C)$  
($\subset \bold C^3$) of  $C$.  If  $u=(u\lr{1}, u\lr{2}, u\lr{3})$  
is in  $\kappa^{-1}\iota(C)$,  we denote by  $(x(u), y(u))$  
the coordinate of the corresponding point on  $C$  by 
  $$
  u\lr{1}=\int_{\infty}^{(x(u), y(u))}\frac{dx}{2y}, \ \ 
  u\lr{2}=\int_{\infty}^{(x(u), y(u))}\frac{xdx}{2y}, \ \ 
  u\lr{3}=\int_{\infty}^{(x(u), y(u))}\frac{x^2dx}{2y} 
  $$
with appropriate choice of a path of the integrals.  
Needless to say, we have  $(x(0, 0, 0), y(0, 0, 0))=\infty$.   

From our standing point of view, 
the following three features stand out on the formulae (1.3) and (1.4).  
Firstly, the sequence of functions of  $u$  
whose values at  $u=u_j$  are displayed in the  $(j+1)$-th row 
of the determinant of (1.3)  is a sequence of the monomials 
of  $x(u)$  and  $y(u)$  displayed according to the order 
of their poles at  $u=0$.  
Secondly, while the right hand sides of (1.3) and (1.4) are 
polynomials of  $x(u)$  and  $y(u)$,  where  $u=u_0$  for (1.4),  
the left hand sides are expressed in terms of theta functions, 
whose domain is properly the universal covering space  $\bold C$  
(of the Jacobian variety) of the elliptic curve.  
Thirdly, the expression of the left hand side of (1.4) states 
the function of the two sides themselves of (1.4) is characterized 
as a hyperelliptic function such that its zeroes are exactly the points 
different from  $\infty$  
whose $n$-plication is just on the standard theta divisor 
in the Jacobian of the curve, and such that its pole is only at  $\infty$.  
In the case of the elliptic curve above, 
the standard theta divisor is just the point at infinity.  

Surprisingly enough, these three features just invent good generalizations 
for hyperelliptic curves.  
More concreatly, our generalization of (1.4) is obtained by replacing 
the sequence of the right hand side by the sequence
  $$
  1,\  x(u),\  x^2(u),\  x^3(u),\  y(u),\  x^4(u),\  yx(u),\  \cdots,
  $$  
where  $u=(u\lr{1}, u\lr{2}, u\lr{3})$  is on  $\kappa^{-1}\iota(C)$, 
of the monomials of  $x(u)$  and  $y(u)$  displayed according to 
the order of their poles at  $u=(0, 0, 0)$  with replacing the derivatives 
with respect to  $u\in\bold C$  by 
those with respect to  $u\lr{1}$  along  $\kappa^{-1}\iota(C)$;  
and the left hand side of (1.4) by  
  $$
  1!2!\cdots (n-1)!\sigma(nu)/\sigma_2(u)^{n^2}, 
  $$  
where  $\sigma(u)=\sigma(u\lr{1}, u\lr{2}, u\lr{3})$  
is a well-tuned Riemann theta series 
and  $\sigma_2(u)=(\partial\sigma/\partial u\lr{2})(u)$.  
Therefore, the hyperelliptic function that is the right hand side of 
the generlarization of (1.4) 
is naturally extended to a function on  $\bold C^3$ 
via theta functions.  
Although the extended function on  $\bold C^3$  is no longer 
a function on the Jacobian, 
it is expressed simply in terms of theta functions and is treated 
really similar way to the elliptic functions.  
The most difficult problem is to find the left hand side of 
the expected generalization of (1.3).  
The answer is remarkably pretty and is 
  $$
  \frac{\sigma(u_0+u_1+\cdots+u_n)\prod_{i<j}\sigma_3(u_i-u_j)}
        {\sigma_2(u_0)^{n+1}\sigma_2(u_1)^{n+1}\cdots\sigma_2(u_n)^{n+1}}, 
  $$
where  $u_j=(u\lr{1}_j,  u\lr{2}_j,  u\lr{3}_j)$  
are variables on  $\kappa^{-1}\iota(C)$  
and  $\sigma_3(u)=(\partial\sigma/\partial u\lr{3})(u)$.  
If we once find this, we can prove the formula, roughly speaking, 
by comparing the divisors of the two sides.  
As same as the formula (1.4) is obtained by 
a limitting process from (1.3), our generalization of (1.4) is 
obtained by similar limitting process from the generalization of (1.3). 

Although this paper is almostly based on \cite{{\bf 13}},  
several critical facts are appeared in comparison with \cite{{\bf 13}}.  
Sections 3 and 4 are devoted to generalize (1.3) and (1.4), respectively.  
We recall in Section 2 the necessary facts for Sections 3 and 4. 

The author started this work by suggestion of S. Matsutani 
concerning the paper \cite{{\bf 13}}.  
After having worked out this paper, 
the author tried to generalize the formula (1.3) further 
to the case of genus larger than three and did not succeed.  
The author hopes that publication of this paper would 
contribute to generalize our formula of type (1.3) or (1.1) 
to cases of genus larger than three in the line of our investigation.  
Matsutani also pointed out that (1.4) can be 
generalized to {\it all} of hyperelliptic curves.  
The reader who is interested in the genaralization of (1.4) 
should be consult with the paper \cite{{\bf 10}}.  

Cantor \cite{{\bf 5}} gave another determiant expression of 
the function that is characterized in the third feature above 
for any hyperelliptic curve.  
The expression of Cantor should be seen as a generalization 
of a formula due to Brioschi (see \cite{{\bf 4}}, p.770, $\ell$.3).  
The Appendix of \cite{{\bf 10}} written by Matsutani 
reveals the connection of our formula, that is Theorem 4.2 below, and 
the determinant expression of \cite{{\bf 5}}.  
So we have three different proofs for the generalization of 
(1.4) in the case of genus three or below.  

There are also various generalizations of (1.1) (or (1.3)) 
in the case of genus two different from our line.  
If the reader is interested in them, 
he should be refered to Introduction of \cite{{\bf 13}}.   

We use the following notations throughout the rest of the paper.  
We denote, as usual, by  $\bold Z$ 
and  $\bold C$ 
the ring of rational integers
and the field of complex numbers, respectively.  
In an expression of the Laurent expansion of a function, 
the symbol  $(d^{\circ}(z_1, z_2, \cdots, z_m)\geq n)$  stands for 
the terms of total degree at least  $n$  with respect to 
the given variables  $z_1$, $z_2$, $\cdots$, $z_m$. 
When the variable or the least total degree are clear from the context, 
we simply denote them by  $(d^{\circ}\geq n)$  or the dots \lq\lq $\cdots$". 

For cross references in this paper, we indicate a formula as (1.2),  
and each of Lemmas, Propositions, Theorems and Remarks also as 1.2.  


\newpage

\heading {\rm 2. {\it The Sigma Function in Genus Three}} \endheading
 
In this Section we summarize the fundamental facts used in Sections 3 and 4.  
Detailed treatment of these facts are given 
in \cite{{\bf 1}}, \cite{{\bf 2}} and \cite{{\bf 3}} 
(see also Section 1 of \cite{{\bf 12}}).  
   
Let 
  $$
  f(x)=\lambda_0     + \lambda_1 x   + \lambda_2 x^2 
     + \lambda_3 x^3 + \lambda_4 x^4 + \lambda_5 x^5  
     + \lambda_6 x^6 + \lambda_7 x^7,  
  $$
where  $\lambda_1$, ... , $\lambda_7$  are fixed complex numbers.  
Assume that the roots of  $f(x)=0$  are different from each other.  
Let  $C$  be a smooth projective model of the hyperelliptic curve 
defined by  $y^2=f(x)$.   Then the genus of  $C$  is  $g$.  
We denote by  $\infty$  the unique point at infinity. 
In this paper we suppose that  $\lambda_7=1$.  The set of forms
  $$
    \omega\lr{1}=\frac{dx}{2y}, 
  \ \omega\lr{2}=\frac{xdx}{2y}, 
  \ \omega\lr{3}=\frac{x^2dx}{2y}
  $$
is a basis of the space of differential forms of first kind.  
We fix generators  $\alpha_1$,  $\alpha_2$,  $\alpha_3$, 
$\beta_1$, $\beta_2$, and  $\beta_3$  of 
the fundamental group of  $C$  such that their intersections 
are  $\alpha_i\cdot\alpha_j=\beta_i\cdot\beta_j=0$,  
$\alpha_i\cdot\beta_j=\delta_{ij}$  for  $i$, $j=1$, $2$, $3$.  
If we set
  $$
  \omega' =\left[\matrix\ 
                           \int_{\alpha_1}\omega\lr{1} &  
                           \int_{\alpha_2}\omega\lr{1} &
                           \int_{\alpha_3}\omega\lr{1} \\
                           \int_{\alpha_1}\omega\lr{2} &
                           \int_{\alpha_2}\omega\lr{2} &
                           \int_{\alpha_3}\omega\lr{2} \\
                           \int_{\alpha_1}\omega\lr{3} &
                           \int_{\alpha_2}\omega\lr{3} &
                           \int_{\alpha_3}\omega\lr{3} \\
            \endmatrix\right], \ \ 
  \omega''=\left[\matrix\ \int_{\beta_1}\omega\lr{1}  &
                           \int_{\beta_2}\omega\lr{1}  & 
                           \int_{\beta_3}\omega\lr{1}  \\
                           \int_{\beta_1}\omega\lr{2}  & 
                           \int_{\beta_2}\omega\lr{2}  &
                           \int_{\beta_3}\omega\lr{2}  \\
                           \int_{\beta_1}\omega\lr{3}  & 
                           \int_{\beta_2}\omega\lr{3}  &
                           \int_{\beta_3}\omega\lr{3}  \\
            \endmatrix\right] 
  $$
the lattice of periods of our Abelian functions appearing below 
is given by 
  $$
  \Lambda 
 =\omega' {\left[\matrix \bold Z \\ \bold Z \\ \bold Z \endmatrix\right]}
 +\omega''{\left[\matrix \bold Z \\ \bold Z \\ \bold Z \endmatrix\right]}
  (\subset \bold C^3). 
  $$
Let  $J$  be the Jacobian variety of the curve  $C$.  
We identify  $J$  with the Picard group  $\text{Pic}^{\circ}(C)$  of 
linear equivalence classes of the divisors of degree 0 of  $C$.  
Let  $\text{Sym}^3(C)$  be the symmetric product of three copies of  $C$.  
Then we have a birational map
  $$
  \align
  \text{Sym}^3(C)&\rightarrow \text{Pic}^{\circ}(C)=J \\
  (P_1, P_2, P_3) &\mapsto \text{the class of}\  P_1+P_2+P_3-3\cdot\infty.
  \endalign
  $$
We may also identify (the  $\bold C$-rational points of)  $J$  
with  $\bold C^3/\Lambda$.  
We denote by  $\kappa$  the canonical 
map  $\bold C^3\rightarrow \bold C^3/\Lambda$  
and by  $\iota$  the embedding of  $C$  into  $J$  given by mapping  $P$  
to  $\text{the class of}\ P-\infty$.  
The image of the triples of the form  $(P_1, P_2, \infty)$,  
by the birational map above, is a theta divisor of  $J$, 
and is denoted by  $\Theta$.  
The image  $\iota(C)$  is obviously contained in  $\Theta$. 
We denote by  $O$  the origin of  $J$.  
Obviously  $\Lambda=\kappa^{-1}(O)=\kappa^{-1}\iota(\infty)$.

\proclaim{\indent\sc Lemma 2.1} As a subvariety of  $J$,  
the divisor  $\Theta$  is singular only at the origin of  $J$. 
\endproclaim

\fp
A proof of this fact is seen, for instance, 
in Lemma 1.7.2(2) of \cite{{\bf 12}}. 

Let  
  $$
  \aligned
  \eta\lr{1}
   &=\frac{( \lambda_3 x   +2\lambda_4 x^2 + 3\lambda_5 x^3
           +4\lambda_6 x^4 +5\lambda_7 x^5)                 dx}{2y}, \\ 
  \eta\lr{2}
   &=\frac{( \lambda_5 x^2 +2\lambda_6 x^3 + 3\lambda_7 x^4)dx}{2y}, \\ 
  \eta\lr{3}
   &=\frac{            x^3 dx}{2y}.
  \endaligned
  $$
Then  $\eta\lr{1}$, $\eta\lr{2}$,  and  $\eta\lr{3}$  are differential forms 
of the second kind without poles except at  $\infty$ 
(see \cite{{\bf 1}, p.195, Ex.i} or \cite{{\bf 2}, p.314}).  
We also introduce matrices 
  $$
  \eta'  =\left[\matrix\ \int_{\alpha_1}\eta\lr{1}  &
                         \int_{\alpha_2}\eta\lr{1}  &
                         \int_{\alpha_3}\eta\lr{1} \\
                         \int_{\alpha_1}\eta\lr{2}  &
                         \int_{\alpha_2}\eta\lr{2}  &
                         \int_{\alpha_3}\eta\lr{2} \\
                         \int_{\alpha_1}\eta\lr{3}  &
                         \int_{\alpha_2}\eta\lr{3}  &
                         \int_{\alpha_3}\eta\lr{3} \\
           \endmatrix\right], \ 
  \eta'' =\left[\matrix\ \int_{\beta_1} \eta\lr{1}  &
                         \int_{\beta_2} \eta\lr{1}  &
                         \int_{\beta_3} \eta\lr{1} \\
                         \int_{\beta_1} \eta\lr{2}  &
                         \int_{\beta_2} \eta\lr{2}  &
                         \int_{\beta_3} \eta\lr{2} \\
                         \int_{\beta_1} \eta\lr{3}  &
                         \int_{\beta_2} \eta\lr{3}  &
                         \int_{\beta_3} \eta\lr{3} \\
           \endmatrix\right].
  $$
The modulus of  $C$  is  $Z:={\omega'}^{-1}\omega''$. 
If we set 
  $$
  \delta'' = \left[\frac 12 \ \ \frac 12 \ \ \frac 12\right], \ \ 
  \delta'  = \left[    0    \ \ \frac 12 \ \    1    \right], \
  $$
then the sigma function attached to  $C$  is defined, as in \cite{{\bf 3}}, by 
  $$
  \aligned
  \sigma&(u) = c\exp(-\frac 12u\eta'{\omega'}^{-1}\tp{u}) \\
  &\cdot\sum_{n\in \bold Z^2}
  \exp[2\pi \sqrt{-1}\{\frac 12\tp{(n+\delta'')}Z(n+\delta'')
                      +\tp{(n+\delta'')}({\omega'}^{-1}\tp{u}+\delta')\}]
  \endaligned
  \tag 2.1
  $$
with a constant  $c$. 
This constant  $c$  is fixed by the following lemma.   

\proclaim{\indent\sc Lemma 2.2}
The Taylor expansion of  $\sigma(u)$  at  $u=(0,0,0)$  is, 
up to a multiplicative constant, of the form 
  $$
  \align
  \sigma(u)
  =&u\lr{1}u\lr{3}-{u\lr{2}}^2 
   -\frac {\lambda_0}3 {u\lr{1}}^4
   -\frac {\lambda_1}3 {u\lr{1}}^3{u\lr{2}}
   -       \lambda_2   {u\lr{1}}^2{u\lr{2}}^2 
   -\frac {\lambda_3}3 {u\lr{1}}  {u\lr{2}}^3 \\
  &-\frac {\lambda_4}3            {u\lr{2}}^4 
   -\frac{2\lambda_2}3 {u\lr{1}}^3           {u\lr{3}} 
   -\frac {\lambda_5}3            {u\lr{2}}^3 u\lr{3}
   -\frac {\lambda_6}2 {u\lr{2}}^2           {u\lr{3}}^2 
   +\frac {\lambda_6}6 {u\lr{1}}             {u\lr{3}}^3 \\
  &-\frac {\lambda_7}3            {u\lr{2}}  {u\lr{3}}^3
   +(d^{\circ}\geq 6), \quad (\lambda_7=1),
  \endalign
  $$
with the coefficient of the term  ${u\lr{3}}^6$  
being  $\frac{\lambda_7}{45}$. 
\endproclaim

\fp
Lemma 2.2 is proved in Proposition 2.1.1(3) of \cite{{\bf 12}} 
by the same argument of \cite{{\bf 12}}, p.96.   
We fix the constant  $c$  in (2.1) such that the expansion is 
exactly of the form in 2.2. 

\proclaim{\indent\sc Lemma 2.3} Let  $\ell$  be an element of  $\Lambda$.  
The function  $u \mapsto \sigma(u)$  on  $\bold C^3$  satisfies 
the translational formula
  $$
  \sigma(u+\ell)=\chi(\ell)\sigma(u)\exp L(u+\ell, \ell), 
  $$
where  $\chi(\ell)=\pm 1$  is independent of  $u$,   
$L(u,\  v)$  is a form which is bilinear over the real field 
and  $\bold C$-linear with respect to the first variable  $u$,  
and  $L(\ell_1,\ \ell_2)$  is  $2\pi\sqrt{-1}$  times an integer 
if  $\ell_1$  and  $\ell_2$  are in  $\Lambda$.  \fp
\endproclaim

\fp The detail of 2.3 is given in  \cite{{\bf 12}}, p.286 
and Lemma 3.1.2 of \cite{{\bf 12}}.   

\proclaim{\indent\sc Lemma 2.4} 
{\rm (1)}  The function  $\sigma(u)$  on  $\bold C^3$  
vanishes if and only if  $u\in \kappa^{-1}(\Theta)$.  \fp 
{\rm (2)}  Suppose that  $v_1$, $v_2$, $v_3$  are three points 
of  $\kappa^{-1}\iota(C)$.  
The function  $u\mapsto \sigma(u-v_1-v_2-v_3)$  is identically zero 
if and only if  $v_1+v_2+v_3$  is contained in  $\kappa^{-1}(C)$.  
If the function is not identically zero, 
it vanishes only at  $u=v_j$  modulo  $\Lambda$  
for  $j=1$, $2$, $3$  of order  $1$  
or of multiple order according as coincidence of 
some of the three points.  \fp
{\rm (3)}  Let  $v$  be a fixed point of  $\kappa^{-1}\iota(C)$.  
There exist two points  $v_1$  and $v_2$  of  $\kappa^{-1}\iota(C)$  
such that the function   $u\mapsto \sigma(u-v-v_1-v_2)$  
on  $\kappa^{-1}\iota(C)$  is not identically zero 
and vanishes at  $u=v$  modulo  $\Lambda$  of order  $1$. 
\endproclaim

\fp{\it Proof.} 
The assertion 2.4(1) and (2) is proved in \cite{{\bf 1}}, pp.252-258, 
for instance. 
The assertion (3) obviously follows from (2).  
\qedright
 
We introduce the functions 
  $$
  \wp_{jk}(u)=-\frac{\partial^2}{\partial u\lr{j}\partial u\lr{k}}
               \log\sigma(u), \ 
  \wp_{jk\cdots r}(u)=\frac{\partial}{\partial u\lr{j}}\wp_{k\cdots r}(u)
  $$
which are defined by Baker.  
Lemma 2.3 shows that these functions are periodic 
with respect to the lattice  $\Lambda$.  
By 2.4(1) we know that the functions  $\wp_{jk}(u)$  
and  $\wp_{jk\ell}(u)$  have its poles along  $\Theta$. 
We also use the notation
  $$
  \sigma_j(u)=\frac{\partial}{\partial u\lr{j}}\sigma(u), \ 
  \sigma_{jk\cdots r}(u)
         =\frac{\partial}{\partial u\lr{j}}\sigma_{k\cdots r}(u).  
  $$

Let  $(u\lr{1}, u\lr{2}, u\lr{3})$  be an arbitrary point in $\bold C^3$.  
Then we can find a set of 
three points  $(x_1, y_1)$,  $(x_2, y_2)$,  and  $(x_3, y_3)$  
on  $C$  such that 
  $$
  \aligned
  u\lr{1}&=\int_{\infty}^{(x_1, y_1)}\omega\lr{1} 
         + \int_{\infty}^{(x_2, y_2)}\omega\lr{1} 
         + \int_{\infty}^{(x_3, y_3)}\omega\lr{1}, \\ 
  u\lr{2}&=\int_{\infty}^{(x_1, y_1)}\omega\lr{2} 
         + \int_{\infty}^{(x_2, y_2)}\omega\lr{2}
         + \int_{\infty}^{(x_3, y_3)}\omega\lr{2}, \\
  u\lr{3}&=\int_{\infty}^{(x_1, y_1)}\omega\lr{3} 
         + \int_{\infty}^{(x_2, y_2)}\omega\lr{3}
         + \int_{\infty}^{(x_3, y_3)}\omega\lr{3} \\
  \endaligned
  \tag 2.2
  $$
with certain choices of the three paths in the integrals. 
If  $(u\lr{1}, u\lr{2}, u\lr{3})$  does not belongs 
to  $\kappa^{-1}(\Theta)$, the set of the three points is uniquely 
determined.  In this situation, one can show the following.

\proclaim{\indent\sc Lemma 2.5}
With the notation above, we have
  $$
  \wp_{13}(u)= x_1x_2x_3, \ \ 
  \wp_{23}(u)=-x_1x_2-x_1x_3-x_3x_1, \ \ 
  \wp_{33}(u)=x_1+x_2+x_3. 
  $$
\endproclaim

\fp
For a proof of this, see \cite{{\bf 2}}, p.377. 
This fact is entirely depends on the choice of forms  $\omega\lr{j}$'s  
and  $\eta\lr{j}$'s.

\proclaim{\indent\sc Lemma 2.6}
If  $u=(u\lr{1}, u\lr{2}, u\lr{3})$  is on $\kappa^{-1}\iota(C)$,  
then we have 
  $$
  u\lr{1}=\frac15 {u\lr{3}}^5 + (d^{\circ}(u\lr{3})\geq 6), \ \ 
  u\lr{2}=\frac13 {u\lr{3}}^3 + (d^{\circ}(u\lr{3})\geq 4).
  $$
\endproclaim

\fp
This is mentioned in \cite{{\bf 12}}, Lemma 2.3.2(2).  
If  $u$  is a point on  $\kappa^{-1}\iota(C)$  
the  $x$-  and  $y$-coordinates of  $\iota^{-1}\kappa(u)$  
will be denoted by  $x(u)$  and  $y(u)$, respectively.  
As is shown in Lemma 2.3.1 of \cite{{\bf 12}}, for instance, 
we see the following.  

\proclaim{\indent\sc Lemma 2.7}
If  $u\in \kappa^{-1}\iota(C)$  then 
  $$
  x(u)=\frac 1{{u\lr{3}}^2}+(d^{\circ}\geq -1),\ \  
  y(u)=-\frac 1{{u\lr{3}}^5}+(d^{\circ}\geq -4). 
  $$
\endproclaim

\proclaim{\indent\sc Lemma 2.8}
{\rm (1)} Let  $u$  be an arbitrary point on  $\kappa^{-1}\iota(C)$.  
Then  $\sigma_2(u)$  is  $0$  if and  only if  $u$  belongs 
to  $\kappa^{-1}(O)$.  \fp
{\rm (2)} The Taylor expansion of the function  $\sigma_2(u)$  on  
$\kappa^{-1}\iota(C)$  at  $u=(0, 0, 0)$  is of the form
  $$
  \sigma_2(u) = -{u\lr{3}}^3 + (d^{\circ}(u\lr{3}) \geq 5).
  $$
\endproclaim

\fp
{\it Proof.}
For (1), assume that  $u\in\kappa^{-1}\iota(C)$  
and  $u\not\in\kappa^{-1}(O)$.  
Then we have 
  $$
  \frac{\sigma_1(u)}{\sigma_2(u)}
  =\frac{\wp_{13}(u)}{\wp_{23}(u)}
  =\frac{x_1x_2x_3}{-x_1x_2-x_2x_3-x_3x_1}\Big|_{x_1=x_2=\infty}
  =-x(u), \ \ 
  \frac{\sigma_3(u)}{\sigma_2(u)}
  =\frac{\wp_{33}(u)}{\wp_{23}(u)}=0.  \ \ 
$$
by using 2.4(1) and 2.5.  
Hence it must be  $\sigma_3(u)=0$ by the second formula.  
If  $\sigma_2(u)=0$  then the first formula yields  $\sigma_1(u)=0$.   
This contradicts to 2.1, 2.4(1) and (2).  
So it must be  $\sigma_2(u)\neq 0$.  
The assertion (2) follows from 2.2 and 2.6.  
\qedright

\proclaim{\indent\sc Lemma 2.9}
Let  $u$  be a point on  $\kappa^{-1}(\Theta)$.  
The function $\sigma_3(u)$  vanishes if and 
only if  $u\in \kappa^{-1}\iota(C)$.  
\endproclaim

\fp{\it Proof.} 
We have already proved in the proof of 2.8 
that if  $u\in \kappa^{-1}\iota(C)$  then  $\sigma_3(u)=0$.  
So we prove the converse.  
Assume that  $u\in\kappa^{-1}(\Theta)$, $u\not\in\kappa^{-1}\iota(C)$,   
and  $u$  corresponds to the pair of points  $(x_1,y_1)$  and  $(x_2,y_2)$.  
Then we have 
  $$
  \frac{\sigma_1(u)}{\sigma_3(u)}
  =\frac{\wp_{13}(u)}{\wp_{33}(u)}=-x_1x_2, \ \ 
  \frac{\sigma_2(u)}{\sigma_3(u)}
  =\frac{\wp_{23}(u)}{\wp_{33}(u)}=-x_1-x_2.
  $$
by using 2.4(1) and 2.5.  If  $\sigma_3(u)=0$, then 
the second formula says that  $\sigma_2(u)=0$, and 
the first one says that  $\sigma_1(u)=0$.   
This contradicts to 2.1 by 2.4(1) and (2).  
So it must be  $\sigma_3(u)\neq 0$.  
\qedright

\proclaim{\indent\sc Lemma 2.10} 
Let  $v$  be a fixed point in  $\kappa^{-1}\iota(C)$  different from 
any point of  $\kappa^{-1}(O)$.  
Then the function
  $$
  u \mapsto \sigma_3(u-v)
  $$
vanishes of order  $2$  at  $u=(0, 0, 0)$.  Precisely, one has 
  $$
  \sigma_3(u-v) = \sigma_2(v){u\lr{3}}^2 + (d^{\circ}(u\lr{3}) \geq 3). 
  $$
\endproclaim

\fp{\it Proof.} 
Since  $u-v$  is on  $\Theta$, we have $\sigma(u-v)=0$. 
If we write  $u$  as  $(x_1, y_1)$  and  $v$  as  $(x_2, y_2)$, 
2.4(1), 2.5 and 2.7 imply that
  $$
  \aligned
  \frac{\sigma_3(u-v)}{\sigma_2(u-v)}
  =&\frac{{\sigma_3}^2-\sigma_{33}\sigma}
        {\sigma_2\sigma_3-\sigma_{23}\sigma}
   (u-v) \\
  =&\frac{\wp_{33}}{\wp_{23}}(u-v) \\
  =&-\frac{x_1+x_2+x_3}
         {x_1x_2 + x_2x_3 + x_3x_1}
         \Big|_{x_3=\infty} \\
  =&-\frac1{x_1+x_2} \\
  =&-\frac1{\left(\frac1 { {u\lr{3}}^2 } + \cdots \right) + x_2} \\
  =&-{u\lr{3}}^2 + \cdots.
  \endaligned
  $$
Since  $\sigma_2(-v)=-\sigma_2(v)$, the desired formula follows.  
\qedright


\proclaim{\indent\sc Lemma 2.11} 
Let  $v$  be a fixed point in  $\kappa^{-1}\iota(C)$  different from 
any points in  $\kappa^{-1}(O)$.  
Then the function
  $$
  u \mapsto \sigma_3(u-v)
  $$
on  $\kappa^{-1}\iota(C)$  has a zero of order  $1$  at  $u=v$.  
\endproclaim

\fp{\it Proof.} 
We denote by  $\frac{du\lr{j}}{dx}$  the derivative of 
the function  $u\mapsto u\lr{j}$  on  $\kappa^{-1}\iota(C)$  by  $x(u)$.  
Since
  $$
   \frac{d(u\lr{j}-v\lr{j})}{d(u\lr{1}-v\lr{1})}
  =\frac{d(u\lr{j}-v\lr{j})}{du\lr{j}}
   \frac{du\lr{j}}{du\lr{1}}
   \frac{du\lr{j}}{d(u\lr{1}-v\lr{1})}
  =\frac{du\lr{j}}{du\lr{1}}
  =\frac{du\lr{j}}{dx}\frac{dx}{du\lr{1}}
  =x^{j-1}(u)
  \tag 2.3
  $$
for  $j=1$  and  $2$,  we see 
  $$
  u\lr{j}-v\lr{j}=x^{j-1}(v)(u\lr{1}-v\lr{1})
           +(d^{\circ}(u\lr{1}-v\lr{1}) \geq 2).
  $$
There exist two points  $v_1$  and  $v_2$  in  $\kappa^{-1}\iota(C)$  
such that  
the fuction  $u\mapsto \sigma(u-v-v_1-v_2)$  on  $\kappa^{-1}\iota(C)$  
is not identically zero and vanishes at  $u=v$  of order  $1$  by 2.4(3).  
Let  $m$  be the vanishing order of the function  $u\mapsto u\lr{1}-v\lr{1}$.  
Then the vanishing orders of  $u\mapsto u\lr{j}-v\lr{j}$  ($j=1$, $2$)  
are equal to or larger than  $m$  by (2.3).  
Furthermore the expansion
  $$
  \multline
  \sigma(u-v-v_1-v_2) \\
  =\sigma_1(-v_1-v_2)(u\lr{1}-v\lr{1})
  +\sigma_2(-v_1-v_2)(u\lr{2}-v\lr{2})
  +\sigma_3(-v_1-v_2)(u\lr{3}-v\lr{3}) \\
  +(d^{\circ}(u\lr{1}-v\lr{1}, u\lr{2}-v\lr{2}, u\lr{3}-v\lr{3})\geq 2)
  \endmultline
  $$
shows that the vanishing order of  $u\mapsto \sigma(u-v-v_1-v_2)$  is 
higher than or equal to  $m$.  Hence  $m$  must be  $1$.  
On the other hand, 2.2 and (2.3) imply that
  $$
  \sigma_3(u-v)
  =(u\lr{1}-v\lr{1}) + (d^{\circ}(u\lr{1}-v\lr{1})\geq 2).
  $$
Thus the statement follows. 
\qedright

\proclaim{\indent\sc Lemma 2.12} 
If  $u$  is a point of  $\kappa^{-1}\iota(C)$, 
then 
  $$
  \frac{\sigma_3(2u)}{\sigma_2(u)^4}=-2y(u). 
  $$
\endproclaim

\fp{\it Proof.} We first prove that 
the left hand side is a function on  $\iota(C)$.  
Since  $[2]^*\Theta$, the pull-back by duplication in  $J$  of  $\Theta$,  
is linearly equivalent to  $4\Theta$ as is shown by Collorary 3 of 
\cite{{\bf 11}}, p.59, and an equality  $[-1]^*\Theta=\Theta$,  
the function  $\sigma(2u)/\sigma(u)^4$  is a function on  $J$.  
For  $u\notin \kappa^{-1}\iota(C)$,  
after multiplying 
  $$
  \frac{\wp_{333}(2u)}{\wp_{33}(2u)\wp_{22}(u)^2}
 =\frac{-2{\sigma_3}^3+3{\sigma_3}^2-\sigma_{333}\sigma^2}
       {{\sigma_3}^2 - \sigma_{33}\sigma} (2u) \cdot
  \frac{\sigma^2}{{\sigma_2}^2-\sigma_{22}\sigma}(u)
  $$
to the function  $\sigma(2u)/\sigma(u)^4$, 
bringing  $u$  close to any point of  $\kappa^{-1}\iota(C)$,  
we obtain the left hand side of the desired formula.  
Here we have used the fact 
that  $u\mapsto \sigma_3(2u)$  does not vanish, 
which follows from 2.9.  
Thus the the function  $\sigma_3(2u)/\sigma_2(u)^4$  is 
a function on  $\iota(C)$, that is 
  $$
  \frac{\sigma_3(2(u+\ell))}{\sigma_2(u+\ell)^4}
  =\frac{\sigma_3(2u)}{\sigma_2(u)^4}
  $$
for  $u\in\kappa^{-1}(C)$.  
Lemma 2.8(1)  states this function has its only pole 
at $u=(0, 0, 0)$  modulo  $\Lambda$.  
Lemma 2.2 and 2.8(2) give 
that its Laurent expansion at  $u=(0, 0, 0)$  is 
  $$
  \frac{ 2\left(\frac15{u\lr{3}}^5\right) 
      - \lambda_7\cdot 2
          \left(\frac13{u\lr{3}}^3\right)
          \left(2u\lr{3}\right)^2
      +   \frac{6\lambda_7}{45}
          \left(2{u\lr{3}}\right)^5
      + \cdots}
       {(- {u\lr{3}}^3 + \cdots )^4}
  = \frac2{{u\lr{3}}^7} + \cdots.
  $$
Here we have used the assumption  $\lambda_7=1$.  
Hence the function must be  $-2y(u)$  by 2.7. 
\qedright

\proclaim{\indent\sc Definition-Proposition 2.13} 
Let  $n$  be a positive integer. If  $u\in\kappa^{-1}\iota(C)$, then 
  $$
  \psi_n(u):=\frac{\sigma(nu)}{\sigma_2(u)^{n^2}}
  $$
is periodic with respect to  $\Lambda$.  
In other words it is a function on  $\iota(C)$. 
\endproclaim

\fp
This is proved by a similar argument of 2.12.  
For details, see Proposition 3.2.2 in \cite{{\bf 12}}, p.396.  
By 2.8(2) the function  $\psi_n(u)$  has its only pole 
at  $u=(0, 0, 0)$  modulo  $\Lambda$. 
Hence it is a polynomial of  $x(u)$  and  $y(u)$.  

$ $

\heading
{\rm 3. {\it A Generalization of the Formula of Frobenius and Stickelberger}} 
\endheading

The following formula is a natural generalization 
of the corresponding formula for Weierstrass' functions  $\sigma(u)$  
and  $\wp(u)$, that is (1.3) for  $n=1$.  

\proclaim{\indent\sc Proposion 3.1} 
If  $u$  and  $v$  are two points in  $\kappa^{-1}\iota(C)$, then
  $$
  \frac{\sigma_3(u + v)\sigma_3(u - v)}
        {{\sigma_2(u)}^2 {\sigma_2(v)}^2}
  =\left|\matrix
      1  &  x(u) \\
      1  &  x(v) 
   \endmatrix\right|.
  $$
\endproclaim

\fp{\it Proof.}
If we regard  $u$  to be a variable on  $\bold C^3$,  
the function 
  $$
  u  \mapsto 
  \frac{\sigma(u+v)\sigma(u-v)}{\sigma(u)^2\sigma(v)^2}
  $$
is periodic with respect to the lattice  $\Lambda$.  
Indeed, the theorem of square (\cite{{\bf 11}}, Coroll. 4 in p.59) yields 
the linear equivalence of  $T^*_{v}\Theta+T^*_{-v}\Theta$  
and  $2\Theta$, where  $T^*_v$  denotes the pull-back of 
the translation by  $v$.  After multiplying 
  $$
  -\frac 12 
  \dfrac{
        \dfrac{\wp_{333}}{\wp_{33}}(u+v)
        \dfrac{\wp_{333}}{\wp_{33}}(u-v)
       }
       {\wp_{22}(u)\wp_{22}(v)}
  $$
to the function above, bringing  $u$  and  $v$  close to points 
on  $\kappa^{-1}\iota(C)$,  
we have the left hand side of the claimed formula because 
of  $\sigma(u\pm v)=\sigma(u)=\sigma(v)=0$  by 2.4(1) (or (2)).  
So the left hand side as a function of  $u$  is periodic 
with respect to  $\Lambda$.  
Now we compare divisors modulo  $\Lambda$  of the two sides.  
The left hand side has its only pole 
at  $u=(0, 0, 0)$  modulo  $\Lambda$  by 2.8(1).  
The two zeroes modulo  $\Lambda$  of the two sides are coincide 
by 2.9 (or 2.11). 
Lemmas 2.8(2) and 2.10 gives its Laurent expansion at  $u=(0, 0, 0)$  
as follows
  $$
  \frac{-\sigma_2(v)({u\lr{3}}^2 + \cdots)\sigma_2(v)({u\lr{3}}^2+\cdots)}
        {(-{u\lr{3}}^3 + \cdots)^2\sigma_2(v)^2}
  =-\frac1{{u\lr{3}}^2}+\cdots.
  $$
The leading term of this coincides with that of the right hand side by 2.7.  
Hence the disired formula holds for all  $v$.  
\qedright

Our generalization of the formula (1.3) in Introduction is the following.  

\proclaim{\indent\sc Theorem 3.2}
Let  $n\geq 2$  be an integer.  
Assume that  $u_0$, $u_1$, $\cdots$,  $u_n$  
belong to  $\kappa^{-1}\iota(C)$.  
Then 
  $$
  \frac{\sigma(u_0+u_1+\cdots+u_n)\prod_{i<j}\sigma_3(u_i-u_j)}
        {\sigma_2(u_0)^{n+1}\sigma_2(u_1)^{n+1}\cdots\sigma_2(u_n)^{n+1}}
  $$
is equal to
  $$
   \left|\matrix
     1  &  x(u_0) & x^2(u_0) & x^3(u_0)  & y(u_0)  & x^4(u_0) &
          yx(u_0) & \cdots   &  yx^{(n-5)/2}(u_0)  & x^{(n+3)/2}(u_0) \\
     1  &  x(u_1) & x^2(u_1) & x^3(u_1)  & y(u_1)  & x^4(u_1) &
          yx(u_1) & \cdots   &  yx^{(n-5)/2}(u_1)  & x^{(n+3)/2}(u_1) \\
 \vdots &  \vdots & \vdots   & \vdots    &\vdots   & \vdots   &
          \vdots  & \ddots   & \vdots              & \vdots           \\
     1  &  x(u_n) & x^2(u_n) & x^3(u_n)  & y(u_n)  & x^4(u_n) & 
          yx(u_n) & \cdots   &  yx^{(n-5)/2}(u_n)  & x^{(n+3)/2}(u_n) \\
   \endmatrix\right|
  $$
or 
  $$
   \left|\matrix
     1  &  x(u_0) & x^2(u_0) & x^3(u_0)  & y(u_0)  & x^4(u_0) &
          yx(u_0) & \cdots   &  x^{(n+2)/2}(u_0)   & yx^{(n-4)/2}(u_0) \\
     1  &  x(u_1) & x^2(u_1) & x^(u_1)^3 & y(u_1)  & x^4(u_1) &
          yx(u_1) & \cdots   &  x^{(n+2)/2}(u_1)   & yx^{(n-4)/2}(u_1) \\
 \vdots &  \vdots & \vdots   & \vdots    & \vdots  & \vdots   & 
          \vdots  & \ddots   & \vdots              & \vdots           \\
     1  &  x(u_n) & x^2(u_n) & x^3(u_n)  & y(u_n)  & x^4(u_n) &
          yx(u_n) & \cdots   &  x^{(n+2)/2}(u_n)   & yx^{(n-4)/2}(u_n) \\
   \endmatrix\right|
  $$
according as  $n$  is odd or even.  Here both of the matrices are 
of size  $(n+1)\times (n+1)$.  
\endproclaim

We prove this Theorem by induction on  $n$.  
First of all we prove the cases of  $n=2$  and  $n=3$.   
We quote these cases as Lemmas 3.3  and  3.4  below.   

\proclaim{\indent\sc Lemma 3.3} 
Assume that  $u$, $u_1$  and  $u_2$  are 
belong to  $\kappa^{-1}\iota(C)$.  Then 
  $$
  \frac{\sigma(u+u_1+u_2)
         \sigma_3(u-u_1)\sigma_3(u-u_2)\sigma_3(u_1-u_2)}
        {\sigma_2(u)^3\sigma_2(u_1)^3\sigma_2(u_2)^3}
  =\left|\matrix
    1 &   x(u)  &   x^2(u)   \\
    1 &  x(u_1) &   x^2(u_1) \\
    1 &  x(u_2) &   x^2(u_2) \\
   \endmatrix\right|.
  $$
\endproclaim

\fp{\it Proof.} 
We suppose that  $u$,  $u_1$,  $u_2$  are any points 
not on  $\kappa^{-1}\iota(C)$.  
Since the sum of pull-backs of
translations  $T^*_{u_1+u_2}\Theta+T^*_{-u_1}\Theta+T^*_{-u_2}\Theta$   
is linearly equivalent to  $3\Theta$  by 
the theorem of square (\cite{{\bf 11}}, Coroll. 4 in p.59),  the function 
  $$
  \frac{\sigma(u+u_1+u_2)
        \sigma(u-u_1)
        \sigma(u-u_2)
        \sigma(u_1-u_2)}
       {\sigma(u)^3\sigma(u_1)^3\sigma(u_2)^3}
  $$
of  $u$  is periodic with respect to the lattice  $\Lambda$. 
As in the proof of 3.1, after multiplying 
  $$
  \frac{
       \dfrac{\wp_{333}}{\wp_{33}}(u-u_1)
       \dfrac{\wp_{333}}{\wp_{33}}(u-u_2)
       \dfrac{\wp_{333}}{\wp_{33}}(u_1-u_2)
       }
       {
       \dfrac{\wp_{222}}{\wp_{22}}(u)
       \dfrac{\wp_{222}}{\wp_{22}}(u_1)
       \dfrac{\wp_{222}}{\wp_{22}}(u_2)
       }
  $$
to the function above, by bringing  $u$,  $u_1$, and  $u_2$  close 
to points on  $\kappa^{-1}\iota(C)$,  
we have the left hand side of the claimed furmula.  
Here we have used the fact that  $\sigma(u-u_1)$, $\sigma(u-u_2)$, 
and  $\sigma(u_1-u_2)$  vanish for  $u$,  $u_1$,  and  $u_2$  on  
$\kappa^{-1}\iota(C)$  by Lemma 2.4(2).  
So the left hand side as a function of  $u$  on  $\kappa^{-1}\iota(C)$  
is periodic with respect to  $\Lambda$.  
Now we regard the both sides to be functions of  $u$  
on  $\kappa^{-1}\iota(C)$.  
We see the left hand side has its only pole at  $u=(0, 0, 0)$  
modulo  $\Lambda$  by 2.8(1), 
and has its zeroes at  $u=\pm u_1$  and  $u=\pm u_2$  modulo  $\Lambda$  
by  2.4(2), 2.9. 
These all zeroes are of order  $1$  by  2.4(2)  and 2.11. 
Its Laurent expansion at  $u=(0, 0, 0)$  is 
given by 2.8(2) and 2.10 as follows:
  $$
  \frac{\sigma_3(u_1+u_2)\sigma_2(u_1)\sigma_2(u_2)\sigma_3(u_1-u_2)}
        {\sigma_2(u_1)^3\sigma_2(u_2)^3}
  (\frac1{{u\lr{3}}^4} + \cdots).
  $$
The right hand side is
  $$
  \left|\matrix 
  1 & x(u_1) \\
  1 & x(u_2)
  \endmatrix\right|\left(\frac1{{u\lr{3}}^4} + \cdots\right).  
  $$
Hence the leading terms of these expansions coincide by 3.1,  
and the sides must be equal. 
\qedright

\proclaim{\indent\sc Lemma 3.4} 
Assume that  $u$, $u_1$, $u_2$  and  $u_3$ 
belong to  $\kappa^{-1}\iota(C)$.  Then 
  $$
  \aligned
  &
  \frac{\sigma(u+u_1+u_2+u_3)
         \sigma_3(u  -u_1)\sigma_3(u  -u_2)\sigma_3(u  -u_3)
         \sigma_3(u_1-u_2)\sigma_3(u_1-u_3)
         \sigma_3(u_2-u_3)}
        {\sigma_2(u)^4\sigma_2(u_1)^4\sigma_2(u_2)^4\sigma_2(u_3)^4} \\
  &=\left|\matrix
    1 &   x(u)  &   x^2(u)   &  x^3(u)  \\
    1 &  x(u_1) &   x^2(u_1) &  x^3(u_1)\\
    1 &  x(u_2) &   x^2(u_2) &  x^3(u_2)\\
    1 &  x(u_3) &   x^2(u_3) &  x^3(u_3)
   \endmatrix\right|.
  \endaligned
  $$
\endproclaim

\fp{\it Proof.} 
We know the left hand side of the claimed formula is, as a function of  $u$,  
a periodic function with respect to  $\Lambda$. 
Its pole is only at  $u=(0, 0, 0)$  modulo  $\Lambda$  and 
is contributed only by the functions 
$\sigma_2(u)^4$, $\sigma_3(u-u_1)$, $\sigma_3(u-u_2)$, $\sigma_3(u-u_3)$.  
By 2.8(2) and 2.10, the order of the pole 
is  $4\times 3 - 3\times 2$, that is  $6$.  
The zeroes of the left hand side are at  $u=-u_1$, $-u_2$, and  $u_3$  
modulo  $\Lambda$  which are coming from  $\sigma(u+u_1+u_2+u_3)$;  
and at  $u=u_1$, $u_2$, $u_3$  which are coming 
from  $\sigma(u-u_1)$, $\sigma(u-u_2)$, $\sigma(u-u_3)$, respectively.  
These  $6$  zeroes are of order  $1$  by 2.11.   
Thus we see that the divisors of two sides coincide.  
The coefficient of leading term of the Laurent expansion of 
the left hand side is 
  $$
  \frac{\sigma(u_1+u_2+u_2)\sigma_2(u_1)\sigma_2(u_2)\sigma_2(u_3)
  \prod_{i<j}\sigma_3(u_i-u_j)}
  {\sigma_2(u_1)^4\sigma_2(u_2)^4\sigma_2(u_3)^4}
  $$
by 2.10.  Such coefficient for the right hand side is
  $$
  \left|\matrix
    1 &  x(u_1) &   x^2(u_1) \\
    1 &  x(u_2) &   x^2(u_2) \\
    1 &  x(u_3) &   x^2(u_3) 
  \endmatrix\right|.
  $$
These two are known to be equal by 3.3 and disired formula is proved.  
\qedright

\fp{\it Proof of Theorem} {\sl 3.2}.  
The best way to explain the general step of the induction 
is probably to demonstrate only the case  $n=4$.  
The case of  $n=4$  is claimed as follows.  
Assume that  $u$, $u_1$, $u_2$, $u_3$, and $u_4$  belong to  $\iota(C)$.  
Then we want to prove the equality 
  $$
  \aligned
  &
  \frac{\sigma(u+u_1+u_2+u_3+u_4)
         \sigma_3(u  -u_1)\sigma_3(u  -u_2)\sigma_3(u  -u_3)\sigma_3(u-u_4)
         \prod_{i<j}\sigma_3(u_i-u_j)}
 {\sigma_2(u)^5\sigma_2(u_1)^5\sigma_2(u_2)^5\sigma_2(u_3)^5\sigma_2(u_4)^5} \\
  &=\left|\matrix
    1 &   x(u)  &   x^2(u)   &  x^3(u)   & y(u)   \\
    1 &  x(u_1) &   x^2(u_1) &  x^3(u_1) & y(u_1) \\
    1 &  x(u_2) &   x^2(u_2) &  x^3(u_2) & y(u_2) \\
    1 &  x(u_3) &   x^2(u_3) &  x^3(u_3) & y(u_3) \\
    1 &  x(u_4) &   x^2(u_4) &  x^3(u_4) & y(u_4)
   \endmatrix\right|.
  \endaligned
  $$
We obviously see that the left hand side of the formula above, 
as a function of  $u$,  is periodic with respect to  $\Lambda$  
by the same argument of 3.1, 3.3 and 3.4,   
and that it has its only pole at  $u=(0, 0, 0)$  modulo  $\Lambda$.  
The order of the pole is  $5\times 3$  coming from  $\sigma_2(u)^5$  
minus  $2\times 4$  coming from  $\sigma_3(u-u_j)$  
for  $j=1$, $2$, $3$, $4$; and that is equal to  $7$.  
We know, by 2.11 that 
there are four obvious zeroes at  $u=u_j$  modulo  $\Lambda$  of order $1$  
coming from  $\sigma(u-u_j)$.  
These are also zeroes of the right hand side.  
Since the right hand side is a polynomial of  $x(u)$  and  $y(u)$,  
it has its only pole at  $u=(0, 0, 0)$  modulo  $\Lambda$.  
Its order is  $7$  coming from the $(1,1)$-entry  $y(u)$.  
So we denote rest zeroes modulo  $\Lambda$  
of the right hand side by  $\alpha$, $\beta$, and $\gamma$.  
Then the theorem of Abel-Jacobi implies 
that  $u_1+u_2+u_3+u_4+ \alpha + \beta + \gamma = (0, 0, 0)$  
modulo  $\Lambda$.  
This means  $\sigma(u+u_1+u_2+u_3+u_4)$  is equal to  
$\sigma(u-\alpha-\beta-\gamma)$  times a trivial theta function.  
Especially these two sigma functions have the same zeroes.  
Since the latter function has obviously zeroes 
at  $u=\alpha$, $\beta$, and  $\gamma$  modulo  $\Lambda$  by 2.4(2), 
the divisors modulo  $\Lambda$  of two sides coincide.  
We can show, as in the proof of 3.3 or 3.4,  
that the coefficients of the leading terms of the two sides 
in their Laurent expansions also coincide 
by using the formula of 3.4.  
The general steps in the induction is done by the same way.  
Thus our proof is completed.  
\qedright

$ $

\centerline{4. {\it Determinant Expression of Generalized Psi-Functions}} 

$ $

In this section we mention a generalization of the formula of (0.2) 
displayed in Introduction.  
Our formula is natural generalization of the formula given in 
Section 3 of \cite{{\bf 13}}.  
Although we can extend this generalization further 
to all of hyperelliptic curves as in \cite{{\bf 10}},  
we give here the case of genus three by limitting process from 3.2.  

The following formula is analogous to 3.1 in \cite{{\bf 13}}.  

\proclaim{\indent\sc Lemma 4.1} Let  $j$  be  $1$, $2$,  or $3$.  
We have 
  $$ 
  \lim_{u\lr{j}\to v\lr{j}}
  \frac{\sigma_3(u-v)}{u\lr{j}-v\lr{j}}=\frac 1{x^{j-1}(v)}.
  $$
\endproclaim

\fp{\it Proof.} Because of 3.1 we have
  $$
  \frac{x(u)-x(v)}{u\lr{j}-v\lr{j}}
  =\frac{\sigma_3(u+v)}{\sigma_2(u)^2\sigma_2(v)^2}
   \cdot
   \frac{\sigma_3(u-v)}{u\lr{j}-v\lr{j}}.
  $$
Now we bring  $u\lr{j}$  close to  $v\lr{j}$.  
Then the limit of the left hand side is  
  $$
  \lim_{u\lr{j}\to v\lr{j}}\frac{x(u)-x(v)}{u\lr{j}-v\lr{j}}
  =\frac{dx}{du\lr{j}}(v).
  $$
This is equal to  $\dfrac{2y}{x^{j-1}}(v)$  by (2.2).  
The assertion follows from 2.12. 
\qedright

Since our proof of the follwing Theorem 
obtained by quite similar argument by using 4.1 
as in the case of genus two (see \cite{{\bf 13}}), 
we leave the proof to the reader.  

\proclaim{\indent\sc Theorem 4.2} 
Let  $n$  be an integer greater than $3$.  
Let  $j$  be any one of  $\{1, 2, 3\}$.  
Assume that  $u$  belongs to  $\kappa^{-1}\iota(C)$.  
Then the following formula for the function  $\psi_n(u)$  
of  {\rm 2.13} holds{\rm :}  
  $$
  \aligned
  &(1!2!\cdots (n-1)!)\psi_n(u) 
  = x^{(j-1)n(n-1)/2}(u)\times \\
  &\left|\matrix 
        x'        & (x^2)'         & (x^3)'     & y'    
  & (x^4)'        & (yx)'          & (x^5)'     & \cdots      \\
       x''        & (x^2)''        & (x^3)''    & y''  
  & (x^4)''       & (yx)''         & (x^5)''    & \cdots      \\
       x'''       & (x^2)'''       & (x^3)'''   & y''' 
  & (x^4)'''      & (yx)'''        & (x^5)'''   & \cdots      \\
     \vdots       & \vdots         & \vdots     & \vdots  
  & \vdots        & \vdots         & \vdots     & \ddots      \\
    x\lr{n-1}     & (x^2)\lr{n-1}  & (x^3)\lr{n-1} & y\lr{n-1} 
  & (x^4)\lr{n-1} & (yx)\lr{n-1} & (x^5)\lr{n-1} & \cdots     \\
  \endmatrix\right|(u).
  \endaligned
  $$
Here the size of the matrix is  $n-1$  by  $n-1$.  
The symbols  ${}'$, ${}''$, $\cdots$, ${}\lr{n-1}$  denote
$\frac{d}{du\lr{j}}$, $\left(\frac{d}{du\lr{j}}\right)^2$, $\cdots$, 
$\left(\frac{d}{du\lr{j}}\right)^{n-1}$, respectively.     
\endproclaim

$ $

{\sl Faculty of Humanities and Social Sciences

Iwate University

Morioka 

020-8550 

Japan} 

$ $

onishi\@iwate-u.ac.jp


\enddocument
\bye